# On the generalized vectorization and its inverse

Vitor Probst Curtarelli

*Abstract*—Although the vectorization operation is known and well-defined, it is only defined for 2-D matrices, and its inverse isn't as well-popularized. This work proposes to generalize the vectorization to higher dimensions, and define mathematically its inverse operation.

*Index Terms*—Vectorization, inverse vectorization, matrices.

## I. Introduction

The vectorization operator [1] is a well-known and defined operator for 2-D matrices, which takes a matrix and turns it into a vector. Its inverse isn't as well spread and aknowledged, although it is can be defined.

This paper proposes to generalize the vectorization operation to higher dimension matrices through different approaches, as well as properly defining the inverse vectorization for such generalizations.

## II. Vectorization

Let $\mathbf{A}$ be a $M \times N$ matrix, such that

$$\mathbf{A} = \begin{bmatrix} a_{0,0} & \cdots & a_{0,N-1} \\ \vdots & \ddots & \vdots \\ a_{M-1,0} & \cdots & a_{M-1,N-1} \end{bmatrix} \quad (1)$$

where the matrix is 0-indexed to facilitate the algebra ahead.

We define $\text{vec}(\cdot)$ as the vectorization operator, given by the vertical stacking of the columns of $\mathbf{A}$. That is, $\mathbf{a} = \text{vec}(\mathbf{A})$ is an $MN \times 1$ vector, given by

$$\mathbf{a} = \text{vec}(\mathbf{A}) = \begin{bmatrix} a_{0,0} \\ \vdots \\ a_{M-1,0} \\ a_{0,1} \\ \vdots \\ a_{M-1,N-2} \\ a_{0,N-1} \\ \vdots \\ a_{M-1,N-1} \end{bmatrix} \quad (2)$$

in which each element of $\mathbf{a}$ is equivalent to one element of $\mathbf{A}$.

## III. Generalized vectorization

Formally, the vectorization is a transformation which is applied to 2-D matrices. However, there are uses where the vectorization of matrices of higher dimension is useful (such as machine learning). In here we will use the loose definition of the vectorization explained previously to more formally define the general vectorization for any dimension.

We will now let $\mathbf{A}$ as a $M_1 \times \cdots M_k$ $k$-dimensional matrix. We also denote $\mathbf{A}^{\mathsf{T}(m,n)}$ as the transpose between the $m$-th and $n$-th dimensions of $\mathbf{A}$.

### A. Block operation

We define the block operation as $\overline{\mathbf{A}} = \mathrm{B}_{\mathbf{T}}(\mathbf{A})$, where $\overline{\mathbf{A}}$ is an $T_1 \times \cdots \times T_k$ block matrix (and $\mathbf{T} = [T_1, \cdots, T_k]$), with each element of $\overline{\mathbf{A}}$ being a block matrix of size $M_1/T_1 \times \cdots \times M_k/T_k$ (assuming that $T_n | M_n \ \forall \ n \in \{1,k\}$). We also define $\mathbf{A} = \mathrm{B}^{-1}(\overline{\mathbf{A}})$ the inverse blocking operation, in which we concatenate the blocks.

### B. Shifting operation

In general, the vectorization transforms a $M \times N$ matrix into a $MN \times 1$ matrix (or a $MN$ vector), "shifting" the last dimension to the one before. Given $\mathbf{A}$ a $k$-D array, if we shift the last dimension to the previous one $k-1$ times, we will achieve a vectorized version of the initial array.

We define the shifting operation as

$$\mathrm{S}_k(\mathbf{A}) = \mathrm{B}^{-1}\left(\mathrm{B}_{1,1,\cdots,M_k}(\mathbf{A})^{\mathsf{T}(k-1,k)}\right) \quad (3)$$

The block operation will produce a $[1 \times 1 \times \cdots \times M_k]$ block matrix. When transposed over its last two dimensions, it becomes a $[1 \times 1 \times \cdots \times M_k \times 1]$ block matrix. The inverse block operation will then concatenate all the blocks, resulting in a $[M_1 \cdots \times (M_{k-1}M_k) \times 1]$ $k$-dimensional matrix, or a $[M_1 \times \cdots \times (M_{k-1}M_k)]$ $(k-1)$-dimensional matrix.

**Example:** Let $\mathbf{A}$ be a $2 \times 2 \times 3$ matrix, given by

$$\mathbf{A} = \left\{ \begin{bmatrix} 1 & 2 & 3 \\ 4 & 5 & 6 \end{bmatrix} \begin{bmatrix} 7 & 8 & 9 \\ 10 & 11 & 12 \end{bmatrix} \right\} \quad (4)$$

Denoting $\mathbf{A}' = \mathrm{B}_{1,1,3}(\mathbf{A})$, then

$$\mathbf{A}' = \left[\left\{\begin{bmatrix}1\\4\end{bmatrix}\ \begin{bmatrix}7\\10\end{bmatrix}\right\}\ \left\{\begin{bmatrix}2\\5\end{bmatrix}\ \begin{bmatrix}8\\11\end{bmatrix}\right\}\ \left\{\begin{bmatrix}3\\6\end{bmatrix}\ \begin{bmatrix}9\\12\end{bmatrix}\right\}\right] \quad (5)$$

where $\mathbf{A}'$ is a $1 \times 1 \times 3$ block matrix, with $2 \times 2 \times 1$ blocks. Denoting $\mathbf{A}'' = \mathbf{A}'^{\mathsf{T}(2,3)}$, then

$$\mathbf{A}'' = \begin{bmatrix}\left\{\begin{bmatrix}1\\4\end{bmatrix}\ \begin{bmatrix}7\\10\end{bmatrix}\right\}\\ \left\{\begin{bmatrix}2\\5\end{bmatrix}\ \begin{bmatrix}8\\11\end{bmatrix}\right\}\\ \left\{\begin{bmatrix}3\\6\end{bmatrix}\ \begin{bmatrix}9\\12\end{bmatrix}\right\}\end{bmatrix} \quad (6)$$

Finally, by applying the inverse blocking, and "simplifying" the $2 \times 6 \times 1$ matrix to a $2 \times 6$, we get

$$\mathrm{S}_3(\mathbf{A}) = \begin{bmatrix}1 & 4 & 2 & 5 & 3 & 6\\ 7 & 10 & 8 & 11 & 9 & 12\end{bmatrix} \quad (7)$$

∎

**Proposition:** *The shifting operation is its own inverse.*

**Proof:** *Let $\mathbf{A}$ be a $[M_1 \times \cdots \times M_k]$ $k$-dimensional matrix, and $\mathbf{A}' = \mathrm{S}_k(\mathbf{A})$. From this, we have that*

$$\begin{aligned}\mathbf{A}'' &= \mathrm{S}_k(\mathbf{A}')\\ &= \mathrm{S}_k(\mathrm{S}_k(\mathbf{A}))\\ &= \mathrm{B}^{-1}\left(\mathrm{B}_{1,\cdots,M_k}\left(\mathrm{B}^{-1}\left(\mathrm{B}_{1,\cdots,M_k}(\mathbf{A})^{\mathsf{T}(k-1,k)}\right)\right)^{\mathsf{T}(k-1,k)}\right)\\ &= \mathrm{B}^{-1}\left(\mathrm{B}_{1,\cdots,M_k}(\mathbf{A})^{\mathsf{T}(k-1,k)\mathsf{T}(k-1,k)}\right)\\ &= \mathrm{B}^{-1}\left(\mathrm{B}_{1,\cdots,M_k}(\mathbf{A})\right)\\ &= \mathbf{A}\end{aligned} \quad (8)$$

*and therefore the shifting operation is its own inverse.*

∎

### C. Generalized block vectorization

We can define the $k$-dimensional vectorization by repeating the shifting operation $k-1$ times, denoted by

$$\mathrm{vec}_k(\mathbf{A}) = \mathrm{S}_k(\mathrm{S}_{k-1}(\cdots \mathrm{S}_2(\mathbf{A}))) \quad (9)$$

We note that, in this case, the repeated shiftings are operating over a decreasing values of $k$ (the first is of $k$, the second of $k-1$, and so forth).

### D. Generalized element-wise vectorization

It is also possible to define a generalized vectorization by calculating the index $m$ in $\mathbf{a}$ referring to the element $(p_1, \cdots, p_k)$ in $\mathbf{A}$.

When calculating the shifting operation, we are creating blocks of decreasingly smaller dimensions of the original matrix. This means that, on the vectorization's output, the values will be grouped by increasingly larger dimensions. Observing the 2-D case, we have that the columns don't "break apart" like the rows do.

We note that we can write $m$ as

$$m = \sum_{l=1}^{k}\left[\left(\prod_{n=0}^{l-1} M_n\right)\left\{\left[m\ /\!/\ \prod_{n=0}^{l-1} M_n\right]\%\ M_l\right\}\right] \quad (10)$$

assuming that $M_0 = 1$, where $/\!/$ and $\%$ are the integer division and remainder operations, respectively. As the term in curly braces is always between 0 and $M_l - 1$, it is replaceable with $p_l$. From this, we then have

$$m = \sum_{l=1}^{k}\left[\prod_{n=0}^{M_n} M_n p_l\right] \quad (11)$$

Therefore, we know the index in $\mathbf{a}$ for each element of $\mathbf{A}$, for matrices of any size.

### E. Row-wise vectorization

Usually, the 2-D vectorization is applied column-wise, however it is also possible to vectorize a matrix row-wise. This is called row vectorization, denoted $\mathrm{rvec}(\mathbf{A})$.

For the 2-D case, we have that $\mathrm{rvec}(\mathbf{A}) = \mathrm{vec}(\mathbf{A}^{\mathsf{T}})$. For larger dimensions, we will have that

$$\mathrm{rvec}_k(\mathbf{A}) = \mathrm{vec}_k\left(\left[\left[\mathbf{A}^{\mathsf{T}(\cdots)}\right]^{\mathsf{T}(2,k-1)}\right]^{\mathsf{T}(1,k)}\right) \quad (12)$$

where we will have $k/\!/2$ transposes between the opposite-index dimensions.

## IV. General inverse vectorization

Given the two definitions for the general vectorization operation, it is possible to define an inverse vectorization operation.

First of all, it is import to note that the vectorization isn't a bijective operation. As an example, it is trivial to see that

$$\mathrm{vec}\left(\begin{bmatrix}a & b\\ c & d\end{bmatrix}\right) = \mathrm{vec}\left(\begin{bmatrix}a & c & b & d\end{bmatrix}\right) \quad (13)$$

However, we seek a way to return to the matrix $\mathbf{A}$, given the vector $\mathbf{a}$. It is possible to show that, given the values of $M_1, \cdots, M_k$, there exists only one matrix $\mathbf{A}$ of this size whose vectorization is $\mathbf{a}$.

## A. Block inverse

Since the shifting operator is invertible (and is its own inverse), then we can invert the vectorization by applying the shifting operation sequentially with increasing index. That is, we define

$$\text{vec}^{-1}_{k;M_1,\cdots,M_k}(\mathbf{a}) = S_2\left(S_3\left(\cdots S_k(\mathbf{a})\right)\right) \\ = S^{-1}_{k;(k-1)}(\mathbf{a}) \quad (14)$$

where $\mathbf{a} = \text{vec}_k(\mathbf{A})$. Logically, we have that $\text{vec}_2(\mathbf{A}) = \text{vec}(\mathbf{A})$, the original vectorization for 2-D matrices.

## B. Element-wise inverse

From eq. (10), we have that

$$p_l = \left[m \mathbin{/\mkern-6mu/} \prod_{n=0}^{l-1} M_n\right] \% M_l \quad (15)$$

and with this, it is easy to calculate the indexes $(p_1, \cdots, p_k)$ of $\mathbf{A}$ for the $m$-th element of $\mathbf{A}$, and with this invert the vectorization process.

## C. Row-wise inverse

Just like we generalized the row vectorization, we can define its inverse, which trivially is

$$\text{rvec}^{-1}_{k;M_1,\cdots,M_k}(\mathbf{a}) = \left(\left(\text{vec}^{-1}_{k;M_k,\cdots,M_1}(\mathbf{a})^{T(\cdots)}\right)^{T(2,k-1)}\right)^{T(1,k)} \quad (16)$$

## D. 2-D Kronecker-product inverse

It is possible to obtain a closed-form formula inverse for the 2-D vectorization operation through the Kronecker product[1]. We let $\mathbf{B}$ be a $MN^2 \times N$ matrix, given by

$$\mathbf{B} = \mathbf{I}_N \otimes \mathbf{a} \quad (17)$$

where $\mathbf{I}_N$ is the $N \times N$ identity matrix, and $\otimes$ denotes the Kronecker product. We also denote $\mathbf{b}_k$ as the $k$-th column of $\mathbf{B}$, and similarly $\mathbf{i}_{N,k}$ is the $k$-th column if $\mathbf{I}_N$. Trivially, we have that $\mathbf{b}_k = \mathbf{i}_{N,k} \otimes \mathbf{a}$, and also $\text{vec}(\mathbf{b}_k) = \mathbf{b}_k$, since it is already a vector.

It is possible to show that $\text{vec}(\mathbf{i}^T_{N,k} \otimes \mathbf{A}) = \mathbf{b}_k$. Using this, and the property $\text{vec}(\mathbf{OPQ}) = (\mathbf{Q}^T \otimes \mathbf{O})\text{vec}(\mathbf{P})$, then

$$\begin{aligned}
\mathbf{A}\mathbf{i}_{N,k} &= \text{vec}(\mathbf{A}\mathbf{i}_{N,k}) \\
&= \text{vec}(\text{vec}(\mathbf{A}\mathbf{I}_N\mathbf{i}_{N,k})) \\
&= \text{vec}((\mathbf{i}^T_{N,k} \otimes \mathbf{A})\text{vec}(\mathbf{I}_N)) \\
&= \text{vec}(\mathbf{I}_M(\mathbf{i}^T_{N,k} \otimes \mathbf{A})\text{vec}(\mathbf{I}_N)) \\
&= \left[\text{vec}(\mathbf{I}_N)^T \otimes \mathbf{I}_M\right]\text{vec}(\mathbf{i}^T_{N,k} \otimes \mathbf{A}) \\
&= \left[\text{vec}(\mathbf{I}_N)^T \otimes \mathbf{I}_M\right]\mathbf{b}_k
\end{aligned} \quad (18)$$

[1]This development was shown by Mårten W. in [2].

By using that matrix multiplication is associative with respect to concatenation, then $\mathbf{A} = \left[\text{vec}(\mathbf{I}_N)^T \otimes \mathbf{I}_M\right]\mathbf{B}$, and with eq. (17) we have

$$\text{vec}^{-1}_{M,N}(\mathbf{a}) = \left[\text{vec}(\mathbf{I}_N) \otimes \mathbf{I}^T_M\right](\mathbf{I}_N \otimes \mathbf{a}) \quad (19)$$